\documentclass[11pt,leqno]{article}
\usepackage{amssymb}
\usepackage{amsmath}
\usepackage{amsthm}
\usepackage{graphicx}
\usepackage{enumerate}
\usepackage[all,2cell]{xy} \UseAllTwocells \SilentMatrices
\reversemarginpar
\numberwithin{equation}{section}
\theoremstyle{definition}

\DeclareMathOperator{\Sp}{Sp}
\DeclareMathOperator{\Mp}{Mp}

\begin{document}
\newcommand{\cb}[1]{\left\{#1\right\}}
\newcommand{\lb}[1]{\left(#1\right)}
\newcommand{\ls}[1]{\left[#1\right]}
\newcommand{\la}[1]{\left|#1\right|}
\newcommand{\bla}[1]{\left|\frac{}{}\!\!#1\right|}
\newcommand{\ld}[1]{\left\|#1\right\|}
\newcommand{\lr}[1]{\left\langle#1\right\rangle}
\newcommand{\lv}[1]{\left|#1\right\rangle}
\newcommand{\rv}[1]{\left\langle#1\right|}
\newcommand{\lrv}[2]{\left\langle#1\right|\left.#2\right\rangle}
\newcommand{\lat}[2]{\left.#1\right|_{#2}}
\newcommand{\dcb}[1]{\cb{\cb{\,#1\,}}}
\newcommand{\dd}[2]{\frac{d{#1}}{d{#2}}}
\newcommand{\ds}[2]{\frac{d^2{#1}}{d{#2}^2}}
\newcommand{\pp}[2]{\frac{\partial{#1}}{\partial{#2}}}
\newcommand{\ps}[2]{\frac{\partial^2{#1}}{\partial{#2}^2}}
\newcommand{\ppp}[3]{\frac{\partial^2{#1}}{\partial{#2}\partial{#3}}}
\newcommand{\ded}[2]{\ensuremath{\frac{\delta{#1}}{\delta{#2}}}}
\newcommand{\tf}{\ensuremath{\tilde{f}}}
\newcommand{\g}{\ensuremath{\mathbf{g}}}
\newcommand{\fg}{\ensuremath{\mathfrak{g}}}
\newcommand{\fh}{\ensuremath{\mathfrak{h}}}
\newcommand{\tp}{\ensuremath{\widetilde{\phi}}}
\newcommand{\hp}{\ensuremath{\widehat{\phi}}}
\newcommand{\tr}{\ensuremath{\widetilde{\rho}}}
\newcommand{\hr}{\ensuremath{\widehat{\rho}}}
\newcommand{\ts}{\ensuremath{\widetilde{\sigma}}}
\newcommand{\hs}{\ensuremath{\widehat{\sigma}}}
\newcommand{\hv}{\ensuremath{\widehat{v}}}
\newcommand{\hvs}{\ensuremath{\widehat{v}^*}}
\newcommand{\hw}{\ensuremath{\widehat{w}}}
\newcommand{\hws}{\ensuremath{\widehat{w}^*}}
\newcommand{\hxi}{\ensuremath{\widehat{\xi}}}
\newcommand{\hx}{\ensuremath{\widehat{x}}}
\newcommand{\hxs}{\ensuremath{\widehat{x}^*}}
\newcommand{\txi}{\ensuremath{\widetilde{\xi}}}
\newcommand{\hy}{\ensuremath{\widehat{y}}}
\newcommand{\hys}{\ensuremath{\widehat{y}^*}}
\newcommand{\sC}{\ensuremath{\mathcal{C}}}
\newcommand{\C}{\ensuremath{\mathbb{C}}}
\newcommand{\E}{\ensuremath{\mathcal{E}}}
\newcommand{\N}{\ensuremath{\mathbb{N}}}
\newcommand{\Q}{\ensuremath{\mathbb{Q}}}
\newcommand{\R}{\ensuremath{\mathbb{R}}}
\renewcommand{\S}{\ensuremath{\mathbb{S}}}
\newcommand{\sS}{\ensuremath{\mathcal{S}}}
\newcommand{\TS}{\ensuremath{TS^3}}
\newcommand{\TpS}{\ensuremath{T^+S^3}}
\newcommand{\TpSn}{\ensuremath{T^+S^3_{n}}}
\newcommand{\Z}{\ensuremath{\mathbb{Z}}}
\newcommand{\contr}{\ensuremath{\lrcorner\,}}
\newcommand{\Det}{\ensuremath{\mbox{Det}}}
\newcommand{\maps}[1]{\ensuremath{\stackrel{#1}{\longrightarrow}}}
\newcommand{\lmaps}[1]{\ensuremath{\stackrel{#1}{\longleftarrow}}}
\newcommand{\opartial}{\ensuremath{\overline{\partial}}}
\newcommand{\fourS}[4]{\ensuremath{\lb{\begin{array}{cc}{#1}&{#2}\\{#3}&{#4}\end{array}}}}
\newcommand{\two}[2]{\ensuremath{\lb{\begin{array}{c}{#1}\\{#2}\end{array}}}}
\newcommand{\twoT}[2]{\ensuremath{\lb{\begin{array}{cc}{#1}&{#2}\end{array}}}}
\newcommand{\three}[3]{\ensuremath{\lb{\begin{array}{c}{#1}\\{#2}\\{#3}\end{array}}}}
\newcommand{\four}[4]{\ensuremath{\lb{\begin{array}{cc}{#1}&{#2}\\{#3}&{#4}\end{array}}}}
\newcommand{\nine}[9]{\ensuremath{\lb{\begin{array}{ccc}{#1}&{#2}&{#3}\\{#4}&{#5}&{#6}\\{#7}&{#8}&{#9}\end{array}}}}

\newtheorem{theorem}{Theorem}[section]
\newtheorem{example}[theorem]{Example}%[section]
\newtheorem{definition}[theorem]{Definition}%[section]
\newtheorem{lemma}[theorem]{Lemma}%[section]

\title{Metaplectic-c Quantized Energy Levels of the Hydrogen Atom}
\author{Jennifer Vaughan\\\texttt{jennifer.vaughan@mail.utoronto.ca}}

\maketitle

\begin{abstract}
This paper calculates the quantized energy levels of the hydrogen atom, using a metaplectic-c prequantization bundle and a definition of a quantized energy level that was introduced by the author in a previous paper \cite{v1}.  The calculation makes use of a computational technique also demonstrated in that paper.  The result is consistent with the standard quantum mechanical prediction.  Unlike other treatments of the hydrogen atom, this approach does not require the construction of the symplectic reduction, but takes place over a regular level set of the energy function.  The Ligon-Schaaf regularization map is used to transform the problem into one of determining the quantized energy levels of a free particle on $\TS$.
\end{abstract}

\section{Introduction}

The quantized energy levels of the hydrogen atom have been calculated for various physical models, using various flavors of geometric quantization.  Notable examples include:
\begin{itemize}
\item Simms \cite{sim3}, who used the observation that the space of orbits corresponding to a fixed negative energy is isomorphic to $S^2\times S^2$, and determined those energies for which the reduced manifold admits a prequantization circle bundle;

\item Sniatycki \cite{sn2}, who looked at the 2-dimensional relativistic Kepler problem and computed a Bohr-Sommerfeld condition for the completely integrable system given by considering the energy and angular momentum functions simultaneously;

\item Duval, Elhadad and Tuynman \cite{det}, who took the phase space to include the spins of the electron and proton, then chose a polarization and determined the Kostant-Souriau quantized operator corresponding to the energy function with fine and hyperfine interaction terms.
\end{itemize}
These examples exist at one of two possible extremes.  On one hand, the quantized energy condition can be evaluated over the symplectic reduction of a particular level set of the energy function, as in \cite{sim3}.  This definition only looks at one energy level at a time, but it requires constructing the symplectic reduction and establishing that the result is a smooth manifold.  On the other hand, the quantized energy levels can be determined from properties of the quantized system as a whole, as in \cite{sn2} or \cite{det}.  These approaches are characterized by requiring a polarization, or the equivalent information -- recall that the Hamiltonian vector fields corresponding to the commuting functions in a completely integrable system generate a real polarization -- and the calculation is not restricted to a single level set of the energy function in question.

In our prior paper \cite{v1}, which builds on work by Robinson \cite{rob1}, we propose an alternative definition of a quantized energy level that acts as a middle ground between these two extremes.  While this definition can be applied to Kostant-Souriau quantization, it was developed for use with metaplectic-c quantization, and that is the context in which we will apply it here.  

Metaplectic-c quantization, due to Robinson and Rawnsley \cite{rr1}, is a quantization procedure that replaces the prequantization circle bundle and metaplectic structure from the Kostant-Souriau recipe by a single object called a metaplectic-c prequantization.  Given a metaplectic-c prequantizable symplectic manifold $(M,\omega)$ and a function $H:M\rightarrow\R$, our quantized energy condition is evaluated over a regular level set of $H$.  It does not require the symplectic reduction, nor does it depend on a polarization.  The objective of this paper is to apply the metaplectic-c quantized energy condition to the hydrogen atom, using the physical model that is equivalent to the Kepler problem.  We will show that the quantized energy levels are in agreement with the quantum mechanical prediction.

In Section \ref{sec:mpc}, we review the fundamentals of metaplectic-c prequantization and the quantized energy condition.  In Section \ref{sec:hatom}, we set up our model of the hydrogen atom and construct a metaplectic-c prequantization for its phase space.  Section \ref{sec:ls} presents the Ligon-Schaaf regularization map, which is a symplectomorphism from the negative-energy domain of the hydrogen atom to an open submanifold of $\TS$.  We show how to relate the energy function for the hydrogen atom to that of a free particle on $\TS$.  Finally, in Section \ref{sec:freeTS}, we determine the quantized energy levels of a free particle on $\TS$, and use these to determine the quantized energy levels for the hydrogen atom.

\section{Metaplectic-c Prequantization and the Quantization Condition}\label{sec:mpc}

In this section, we briefly outline the key features of metaplectic-c prequantization and our quantized energy condition.  A more detailed overview of similar material appears in \cite{v1}, and proofs of the properties of the quantized energy condition are given there.  For the original presentations of metaplectic-c quantization and the constructions of the associated bundles used in the quantized energy condition, see \cite{rr1} and \cite{rob1}.  

Section \ref{subsec:stdvec} establishes standard notation and conventions that will be used throughout the rest of the paper.  Readers who are familiar with the background should nevertheless refer to this section before proceeding to the calculations.

\subsection{The Metaplectic-c Group}\label{subsec:vtspc}

Let $(V,\Omega)$ be a $2n$-dimensional symplectic vector space, and assume that it is equipped with a compatible complex structure $J$.  The metaplectic group $\Mp(V)$ is the unique connected double cover of the symplectic group $\Sp(V)$, and the metaplectic-c group $\Mp^c(V)$ is defined to be $$\Mp^c(V)=\Mp(V)\times_{\Z_2}U(1).$$  

By construction, $\Mp^c(V)$ contains $U(1)$ and $\Mp(V)$ as subgroups.  The inclusion of each subgroup into $\Mp^c(V)$ yields a short exact sequence and a group homomorphism on $\Mp^c(V)$.  One such sequence is $$1\rightarrow U(1)\rightarrow\Mp^c(V)\maps{\sigma}\Sp(V)\rightarrow 1,$$ where the group homomorphism $\sigma$ is called the \emph{projection map}, and its restriction to $\Mp(V)\subset\Mp^c(V)$ is the double covering map.  The other is $$1\rightarrow\Mp(V)\rightarrow\Mp^c(V)\maps{\eta}U(1)\rightarrow 1,$$ where the group homomorphism $\eta$ is called the \emph{determinant map}, and it has the property that for any $\lambda\in U(1)\subset\Mp^c(V)$, $\eta(\lambda)=\lambda^2$.  The Lie algebra $\mathfrak{mp}^c(V)$ is identified with $\mathfrak{sp}(V)\oplus\mathfrak{u}(1)$ under the map $\sigma_*\oplus\frac{1}{2}\eta_*$.

For all $g\in\Sp(V)$, define $$C_g=\frac{1}{2}(g-JgJ).$$  Then $C_g$ is a complex automorphism of $V$.  Using this definition, we construct an embedding of $\Mp^c(V)$ into $\Sp(V)\times\C\setminus\cb{0}$ in the following way.  Given any $a\in\Mp^c(V)$, let $g=\sigma(a)\in\Sp(V)$, and let $\mu\in\C\setminus\cb{0}$ be such that $\eta(a)=\mu^2\Det_\C C_g$.  To remove the ambiguity from the definition of $\mu$, assume that if $a=I$, then $\mu=1$.  The desired embedding is given by the map $a\mapsto(g,\mu)$.  The properties of this map are established in \cite{rr1}; following their terminology, we refer to the pair $(g,\mu)$ as the \emph{parameters} of $a$.  Note that if $a\in\Mp(V)$, then its parameters $(g,\mu)$ satisfy $\mu^2\Det_\C C_g=1$.

Fix a subspace $W\subset V$ with codimension $1$.  Its symplectic orthogonal $W^\perp$ is a one-dimensional subspace of $W$, and the quotient space $W/W^\perp$ acquires a symplectic structure from the symplectic form $\Omega$.  This new symplectic vector space has symplectic group $\Sp(W/W^\perp)$ and metaplectic-c group $\Mp^c(W/W^\perp)$.  The complex structure $J$ on $V$ induces a complex structure on $W/W^\perp$ in such a way that 
\begin{equation}\label{eq:JW}
W/W^\perp\cong(W^\perp\oplus JW^\perp)^\perp,
\end{equation}
where the isomorphism respects both the symplectic structures and the complex structures.

Let $\Sp(V;W)\subset\Sp(V)$ be the subgroup whose elements preserve the subspace $W$.  Then there is a natural group homomorphism $\Sp(V;W)\maps{\nu}\Sp(W/W^\perp)$.  Now let $\Mp^c(V;W)=\sigma^{-1}(\Sp(V;W))\subset\Mp^c(V)$.  As shown in \cite{rob1}, there is a group homomorphism $\Mp^c(V;W)\maps{\hat{\nu}}\Mp^c(W/W^\perp)$ such that the following commutes.
\begin{center}
$\xymatrix{
\Mp^c(V)\supset\Mp^c(V;W) \ar[r]^(0.6){\hat{\nu}} \ar[d]^\sigma & \Mp^c(W/W^\perp) \ar[d]^\sigma\\
\Sp(V)\supset\Sp(V;W) \ar[r]^(0.6){\nu} & \Sp(W/W^\perp)
}$
\end{center}
On the level of Lie algebras, $\eta_*\circ\hat{\nu}_*=\eta_*$.

\subsection{Metaplectic-c Prequantization of Symplectic Manifolds}

Assume that $(V,\Omega)$ has been fixed as in the previous section.  Let $(M,\omega)$ be a $2n$-dimensional symplectic manifold.  We view the symplectic frame bundle $\Sp(M,\omega)\maps{\rho}M$ as a right principal $\Sp(V)$ bundle whose fibers take the form $$\Sp(M,\omega)_m=\cb{b:V\rightarrow T_mM:b\mbox{ is a symplectic isomorphism}},\ \ \forall m\in M.$$  The group action is by precomposition.

\begin{definition}
A \emph{metaplectic-c prequantization} for $(M,\omega)$ is the triple $(P,\Sigma,\gamma)$, where:
\begin{enumerate}[(1)]
\item $P\maps{\Pi}M$ is a right principle $\Mp^c(V)$ bundle over $M$;
\item $\Sigma:P\rightarrow\Sp(M,\omega)$ is a map such that $\rho\circ\Sigma=\Pi$ and $\Sigma(q\cdot a)=\Sigma(q)\cdot\sigma(a)$ for all $q\in P$ and $a\in\Mp^c(V)$; and
\item $\gamma$ is a $\mathfrak{u}(1)$-valued one-form on $P$ such that
\begin{enumerate}[(a)]
\item $\gamma$ is invariant under the $\Mp^c(V)$ action;

\item if $\partial_\alpha$ is the vector field on $P$ generated by $\alpha\in\mathfrak{mp}^c(V)$, then $\gamma(\partial_\alpha)=\frac{1}{2}\eta_*\alpha$; and

\item $d\gamma=\frac{1}{i\hbar}\Pi^*\omega$.
\end{enumerate}
\end{enumerate}
\end{definition}
Note that $(P,\gamma)$ is a principal circle bundle with connection one-form over $\Sp(M,\omega)$, with projection map $\Sigma$.  We will view $P$ as a bundle over $M$ or a bundle over $\Sp(M,\omega)$ as the circumstance demands.

Two metaplectic-c prequantizations for $(M,\omega)$ are considered equivalent if there is a diffeomorphism between them that preserves the one-forms and commutes with the respective maps to $\Sp(M,\omega)$.  If $(M,\omega)$ is metaplectic-c prequantizable, then the set of equivalence classes of metaplectic-c prequantizations for $(M,\omega)$ is in one-to-one correspondence with the locally constant cohomology group $H^1(M,U(1))$.  In particular, if $H^1(M,U(1))$ is trivial, then the metaplectic-c prequantization of $(M,\omega)$ is unique up to isomorphism.

\subsection{Metaplectic-c Quantized Energy Condition}\label{subsec:mpcquant}

Assume that $(M,\omega)$ admits a metaplectic-c prequantization $(P,\Sigma,\gamma)$.  Let $H:M\rightarrow\R$ be a smooth function.  We denote its Hamiltonian vector field by $\xi_H$, and adopt the convention that $\xi_H\contr\omega=dH$.  Let $\xi_H$ have flow $\phi^t$ on $M$.  There is a lift of $\phi^t$ to $\Sp(M,\omega)$ given by $$\tp^t(b)=\phi^t_*\circ b,\ \ \forall b\in\Sp(M,\omega).$$  Let the corresponding vector field on $\Sp(M,\omega)$ be $\txi_H$.  Now let $\hxi_H$ be the lift of $\txi_H$ to $P$ such that $\gamma(\hxi_H)=0$, and let $\hp^t$ be the flow of $\hxi_H$.  

Given a regular value $E$ of $H$, the level set $S=H^{-1}(E)$ is an embedded submanifold of $M$ with codimension $1$.  Assume that a model subspace $W\subset V$ of codimension $1$ has been fixed.  The groups defined in Section \ref{subsec:vtspc} will be the structure groups for the bundles that we now define over $S$.

Let $TS^\perp$ denote the null foliation, which is the bundle over $S$ given by $$T_sS^\perp=\cb{\zeta\in T_sM:\omega(\zeta,T_sS)=0},\ \ \forall s\in S.$$  At each $s\in S$, $T_sS^\perp$ is a one-dimensional subspace of $T_sS$, spanned by $\xi_H(s)$.  The quotient bundle $TS/TS^\perp$ is a symplectic vector bundle over $S$, and we model its symplectic frame bundle $\Sp(TS/TS^\perp)$ on $W/W^\perp$:  $$\Sp(TS/TS^\perp)_s=\cb{b':W/W^\perp\rightarrow T_sS/T_sS^\perp:b'\mbox{ is a symplectic isomorphism}},\ \ \forall s\in S.$$  Let $\Sp(M,\omega;S)\subset\Sp(M,\omega)|_{S}$ be the subbundle defined by  $$\Sp(M,\omega;S)_s=\cb{b\in\Sp(M,\omega):bW=T_sS},\ \ \forall s\in S.$$  This is a principal $\Sp(V;W)$ bundle over $S$.  The symplectic frame bundle $\Sp(TS/TS^\perp)$ is naturally identified with the bundle associated to $\Sp(M,\omega;S)$ by the group homomorphism $\Sp(V;W)\maps{\nu}\Sp(W/W^\perp)$.  

On the level of metaplectic-c bundles, let $P^S$ be the principal $\Mp^c(V;W)$ bundle over $S$ obtained by restricting $P$ to $\Sp(M,\omega;S)$.  Let $P_S$ be the bundle associated to $P^S$ by the group homomorphism $\Mp^c(V;W)\maps{\hat{\nu}}\Mp^c(W/W^\perp)$, and let the projection map to $\Sp(TS/TS^\perp)$ be $P_S\maps{\Sigma_S}\Sp(TS/TS^\perp)$.  The connection one-form $\gamma$ on $P$ pulls back to $\gamma^S$ on $P^S$, and induces a connection one-form $\gamma_S$ on $P_S$.  The following diagram shows the various bundles over $S$, together with their structure groups.
$$
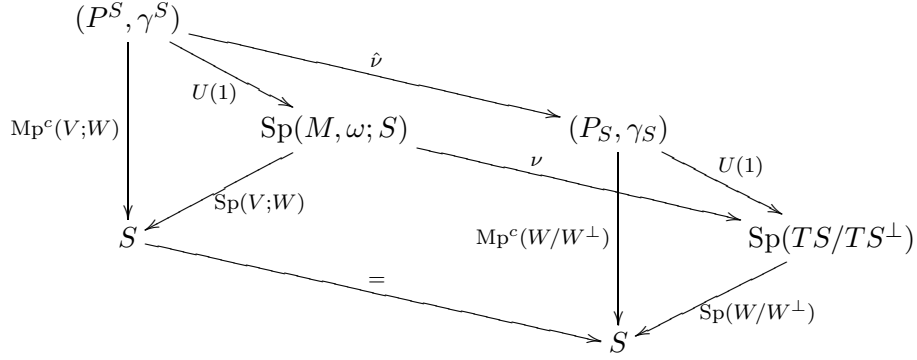
\xymatrix{
(P^S,\gamma^S) \ar[dd]_{\Mp^c(V;W)} \ar[dr]_{U(1)} \ar[drrr]^{\hat{\nu}} \\
& \Sp(M,\omega;S) \ar[dl]^{\Sp(V;W)} \ar[drrr]^(0.4){\nu} & & (P_S,\gamma_S) \ar[dd]_{\Mp^c(W/W^\perp)} \ar[dr]^{U(1)} \\
S \ar[drrr]^{=} & & & & \Sp(TS/TS^\perp) \ar[dl]^{\Sp(W/W^\perp)} \\
& & & S
}$$

It is clear that the flow $\phi^t$ of $\xi_H$ preserves the manifold $S$.  Further, the lifted flows $\tp^t$ and $\hp^t$ preserve the manifolds $\Sp(M,\omega;S)$ and $P^S$, respectively, and they induce flows on $\Sp(TS/TS^\perp)$ and $P_S$.  We also denote these induced flows by $\tp^t$ and $\hp^t$, and let the corresponding vector fields on $\Sp(TS/TS^\perp)$ and $P_S$ be $\txi_H$ and $\hxi_H$.

Having established all of these constructions, we can now state the definition of a quantized energy level for the system $(M,\omega,H)$.  This definition first appeared in \cite{v1}, which contains a detailed examination of its properties.
\begin{definition}\label{def:quantE}
The regular value $E$ of $H$ is a \textbf{quantized energy level} of $(M,\omega,H)$ if the connection one-form $\gamma_S$ on $P_S$ has trivial holonomy over all closed orbits of $\tp^t$ on $\Sp(TS/TS^\perp)$.
\end{definition}

The following property of Definition \ref{def:quantE} was proved in \cite{v1}.
\begin{theorem}
Suppose $H_1,H_2:M\rightarrow\R$ are smooth functions, and suppose there are $E_1,E_2\in\R$ such that $E_j$ is a regular value of $H_j$ for $j=1,2$ and $H_1^{-1}(E_1)=H_2^{-1}(E_2)$.  Then $E_1$ is a quantized energy level for $(M,\omega,H_1)$ if and only if $E_2$ is a quantized energy level for $(M,\omega,H_2)$.
\end{theorem}
In other words, the quantized energy condition depends on the geometry of the level set $S$ and not on a particular choice of $H$, a property that we have termed dynamical invariance.  As a special case of this theorem, suppose $H_1,H_2:M\rightarrow\R$ are smooth functions such that $H_2=f\circ H_1$ for some diffeomorphism $f:\mbox{range}(H_1)\rightarrow\mbox{range}(H_2)$.  Then $E$ is a quantized energy level for $H_1$ if and only if $f(E)$ is a quantized energy level for $H_2$.  We will make use of this fact in Section \ref{subsec:rescale}.

\subsection{Choices for Subsequent Calculations}\label{subsec:stdvec}

In the sections that follow, we will need model symplectic vector spaces of several different dimensions.  Let us fix some standard definitions and notation.  

For any $n\in\N$, let $(V_n,\Omega_n)$ be a $2n$-dimensional symplectic vector space.  Let  $(\hv_1,\ldots,\hv_n,$ $\hw_1,\ldots,\hw_n)$ be a symplectic basis for $V_n$, and write all elements of $V_n$ as ordered $2n$-tuples with respect to this basis.  The symplectic form can be written in terms of the dual basis as $$\Omega_n=\sum_{j=1}^n\hvs_j\wedge\hws_j.$$  Assume that each real vector $(a_1,\ldots,a_n,b_1,\ldots,b_n)\in V_n$ is identified with the complex vector $(b_1+ia_1\ldots,b_n+ia_n)\in\C^n$.  The resulting complex structure $J$ on $V_n$ is written in matrix form as $J=\four{0}{I}{-I}{0}$, where $I$ is the $n\times n$ identity matrix.  

When we require a subspace of $V_n$ of codimension $1$, we choose  $$W_n=\mbox{span}\cb{\hv_1,\ldots,\hv_n,\hw_1,\ldots,\hw_{n-1}}.$$  Then $$W_n^\perp=\mbox{span}\cb{\hv_n}\ \ \mbox{and}\ \ W_n/W_n^\perp=\mbox{span}\cb{[\hv_1],\ldots,[\hv_{n-1}],[\hw_1],\ldots,[\hw_{n-1}]}.$$  Using equation \ref{eq:JW}, it is immediate that $W_n/W_n^\perp$ is isomorphic to $V_{n-1}$ as a symplectic vector space and a complex vector space.  The commutative diagram from Section \ref{subsec:vtspc} containing the group homomorphisms $\nu$ and $\hat{\nu}$ can be rewritten as follows. 
$$\xymatrix{
\Mp^c(V_n)\supset\Mp^c(V_n;W_n) \ar[r]^(0.65){\hat{\nu}} \ar[d]^\sigma & \Mp^c(V_{n-1}) \ar[d]^{\sigma} \\
\Sp(V_n)\supset\Sp(V_n;W_n) \ar[r]^(0.65){\nu} \ar[r] & \Sp(V_{n-1})
}$$

\section{The Hydrogen Atom}\label{sec:hatom}

\subsection{Setup}\label{subsec:setupH}

Let $\dot{\R}^3$ represent $\R^3$ with the origin removed, and let $M=T\dot{\R}^3=\dot{\R}^3\times\R^3$.  We use Cartesian coordinates $q=(q_1,q_2,q_3)$ on $\dot{\R}^3$ and $p=(p_1,p_2,p_3)$ on $\R^3$.  Equip $M$ with the symplectic form $\omega=\sum_{j=1}^3dq_j\wedge dp_j$.

We consider the model of the hydrogen atom that is equivalent to the Kepler problem.  Assume that a proton is fixed at the origin in $\R^3$, and that an electron of mass $m_e$ interacts with it via the electrostatic force, which obeys an inverse-square law with constant of proportionality $k$.  Then $M$ is the phase space for the motion of the electron, where $q$ and $p$ represent its position and momentum, respectively.  The Hamiltonian energy function is $$H=\frac{1}{2m_e}|p|^2-\frac{k}{|q|},$$ and the corresponding Hamiltonian vector field is $$\xi_H=\sum_{j=1}^3\lb{\frac{1}{m_e}p_j\pp{}{q_j}-\frac{k}{|q|^3}q_j\pp{}{p_j}}.$$  Let $\xi_H$ have flow $\phi^t$ on $M$.

The solutions to the Kepler problem are well known \cite{a1,cb}.  The angular momentum vector $L=q\times p$ is a constant of the motion, as is the eccentricity vector $e=\frac{1}{m_ek}p\times L-\frac{q}{|q|}$.  In position space, the orbits corresponding to a given energy $E\in\R$ are conic sections with eccentricity 
\begin{equation}\label{eq:ecc}
|e|=\sqrt{1+\frac{2E|L|^2}{m_ek^2}},
\end{equation}
having the origin as a focus.

In particular, suppose $E<0$.  Then the value of $|L|$ lies in the interval $\ls{0,\sqrt{-\frac{m_ek^2}{2E}}}$.  If $|L|>0$, then the orbit is an ellipse, with $|L|=\sqrt{-\frac{m_ek^2}{2E}}$ being the special case of a circle.  All elliptical orbits with energy $E$ have period $\frac{2\pi}{\Lambda}$, where $\Lambda=\sqrt{-\frac{8E^3}{m_ek^2}}$.  If $|L|=0$, however, then $e=1$ and the orbit is a line segment.  Physically, this represents the case where the electron begins from rest and collapses on a straight-line trajectory into the proton.  Such motion is not periodic, which implies that the level set $H^{-1}(E)$ contains orbits of $\xi_H$ that are not closed.  Further, the collapse occurs in finite time, meaning that the vector field $\xi_H$ is not complete.

The objective of this paper is to determine the quantized energy levels for $(M,\omega,H)$.  In the next section, we construct a metaplectic-c prequantization for $(M,\omega)$, and formulate our approach for performing the quantized energy calculation.

\subsection{Metaplectic-c Prequantization of $(M,\omega)$}\label{subsec:mpcP}

We choose the model symplectic vector space $V_3$, as described in Section \ref{subsec:stdvec}.  The tangent bundle $TM$ can be identified with $M\times V_3$ with respect to the global trivialization $$\hv_j\mapsto\lat{\pp{}{q_j}}{m},\ \ \hw_j\mapsto\lat{\pp{}{p_j}}{m},\ \ \forall m\in M,\ j=1,2,3.$$ This yields an identification of the symplectic frame bundle $\Sp(M,\omega)$ with  $M\times\Sp(V_3)$. 

Let $P=M\times\Mp^c(V_3)$, with bundle projection map $P\maps{\Pi}M$.  Define the map $\Sigma:P\rightarrow\Sp(M,\omega)$ by $$\Sigma(m,a)=(m,\sigma(a)),\ \ \forall m\in M,\forall a\in\Mp^c(V_3),$$ where the right-hand side is written with respect to the global trivialization above.  Let 
\begin{equation}\label{eq:beta}
\beta=\sum_{j=1}^3\lb{q_jdp_j+d(q_jp_j)}
\end{equation}
on $M$, so that $d\beta=\omega$.  The reason for this choice of $\beta$ will be made clear in Section \ref{subsec:lsts}.  Let $\gamma$ be the $\mathfrak{u}(1)$-valued one-form on $P$ given by $$\gamma=\frac{1}{i\hbar}\Pi^*\beta+\frac{1}{2}\eta_*\vartheta_0,$$ where $\vartheta_0$ is the trivial connection on the product bundle.  Then $(P,\Sigma,\gamma)$ is a metaplectic-c prequantization for $(M,\omega)$, and it is unique up to isomorphism.

The quantized energy levels of $(M,\omega,H)$ are those regular values $E$ of $H$ such that the holonomy of $\gamma_S$ is trivial over all closed orbits of $\txi_H$ on $\Sp(TS/TS^\perp)$, where $S=H^{-1}(E)$.  If $E\geq 0$, then the quantization condition can be evaluated immediately.  From equation \ref{eq:ecc} for the eccentricity of the orbit, we see that if $E=0$, then the orbits are parabolas in position space, and if $E>0$, then they are hyperbolas.  In these cases, $\xi_H$ has no closed orbits in $S$, which implies that $\txi_H$ cannot have any closed orbits in $\Sp(TS/TS^\perp)$.  Therefore the holonomy condition is satisfied vacuously, and all nonnegative energy levels are quantized energy levels.  This is consistent with the physical prediction from quantum mechanics:  a particle that is not spatially confined has a continuous energy spectrum.

It remains to consider the orbits corresponding to negative energy.  Let $$N=\cb{m\in M:H(m)<0}.$$ Then $N$ is an open, simply connected submanifold of $M$.  By restriction, the symplectic form on $M$ induces one on $N$, and the metaplectic-c prequantization for $M$ induces one for $N$.  We use the same symbols to denote the restricted objects:  $(N,\omega)$ is a symplectic manifold, and $(P,\Sigma,\gamma)$ is its metaplectic-c prequantization.  Since $N$ is simply connected, $(P,\Sigma,\gamma)$ is unique up to isomorphism.

Let $E<0$ be fixed, and let $S=H^{-1}(E)\subset N$.  Through the process described in Section \ref{subsec:mpcquant}, we obtain the three-level structures $$(P^S,\gamma^S)\rightarrow\Sp(N,\omega;S)\rightarrow S$$ and $$(P_S,\gamma_S)\rightarrow\Sp(TS/TS^\perp)\rightarrow S.$$  Lift $\xi_H$ to $\txi_H$ on $\Sp(N,\omega;S)$, and let it induce a vector field on $\Sp(TS/TS^\perp)$.  

To evaluate the quantization condition using these bundles, we would have to determine the closed orbits of $\txi_H$ on $\Sp(TS/TS^\perp)$, then lift $\txi_H$ to $\hxi_H$ on $P_S$, horizontally with respect to $\gamma_S$, and ensure that every lift of a closed orbit in $\Sp(TS/TS^\perp)$ is closed in $P_S$.  However, this procedure is computationally prohibitive in all but the special case of the circular orbit.  Instead, we will use Ligon-Schaaf regularization to transform the quantized energy calculation on $(N,\omega)$ into one on an open submanifold of $\TS$.  This is the subject of Section \ref{sec:ls}.

\section{Ligon-Schaaf Regularization}\label{sec:ls}

Let $\TpS$ represent the result of removing the zero section from $\TS$.  In Section \ref{subsec:setupH}, we noted that the vector field $\xi_H$ is not complete.  The Ligon-Schaaf map is a symplectomorphism from $(N,\omega)$ to an open submanifold of $\TpS$, having the property that $\xi_H$ is mapped to a vector field that is complete on $\TpS$.  This map was presented in \cite{ls}, and an in-depth discussion of it can be found in \cite{cb}.

In this section, we will state the Ligon-Schaaf map and list its relevant properties.  Then we will construct a metaplectic-c prequantization for $\TS$, and show how to lift the Ligon-Schaaf map to the level of symplectic frame bundles, and of metaplectic-c prequantizations.  Using the lifted maps, we will be able to relate the quantized energy levels of the hydrogen atom to those of a free particle on $\TS$.

Our conventions and notation largely follow those in \cite{cb}.  We state results without proof; much more detail can be found in \cite{cb} and works cited.

\subsection{$\TS$ and the Ligon-Schaaf Map}\label{subsec:lsts}

Consider $T\R^4=\R^4\times\R^4$ with Cartesian coordinates $(x,y)$, where $x=(x_1,x_2,x_3,x_4)$ and $y=(y_1,y_2,y_3,y_4)$.  Let $T\R^4$ have symplectic structure $$\omega_4=\sum_{j=1}^4 dx_j\wedge dy_j,$$ and let 
\begin{equation}\label{eq:beta4}
\beta_4=-\sum_{j=1}^4y_jdx_j
\end{equation}
on $T\R^4$, so that $d\beta_4=\omega_4$.  The submanifold $\TS$ is given by $$\TS=\cb{(x,y)\in T\R^4:|x|^2=1,x\cdot y=0}\subset T\R^4,$$ where we take the usual Euclidean inner product on $\R^4$.  From now on, we abbreviate $T\R^4$ by $Z$.

By treating $(\TS,\omega_3)$ as a level set of the constraint functions $c_1(x,y)=\frac{1}{2}(|x|^2-1)$ and $c_2(x,y)=x\cdot y$ on  $(Z,\omega_4)$, it can be shown that the restriction of $\omega_4$ to $\TS$ yields a symplectic form on $\TS$.  Let $(\TS,\omega_3)$ be the resulting symplectic manifold.  Let $\beta_3$ be the restriction of $\beta_4$ to $\TS$, so that $d\beta_3=\omega_3$.

Let $\TpS$ represent $\TS$ with the zero section removed:  $$\TpS=\cb{(x,y)\in\TS:|y|>0}.$$  Define the map $D:\TpS\rightarrow\R$ by $$D(x,y)=-\frac{m_ek^2}{2|y|^2},\ \ \forall (x,y)\in\TpS.$$  %By using the map $D(x,y)=-\frac{m_ek^2}{2|y|^2}-\frac{1}{2}(x\cdot y)^2+\frac{m_ek^2}{2|y|^2}(|x|^2-1)$ on $T\R^4$, 
As shown in \cite{cb}, the Hamiltonian vector field for $D$ on $\TpS$ is  $$\xi_D=\sum_{j=1}^4\lb{\frac{m_ek^2}{|y|^4}y_j\pp{}{x_j}-\frac{m_ek^2}{|y|^2}x_j\pp{}{y_j}},$$ and its flow is $$\psi^t_D(x,y)=\four{\cos\frac{m_ek^2}{|y|^3}t}{\frac{1}{|y|}\sin\frac{m_ek^2}{|y|^3}t}{-|y|\sin\frac{m_ek^2}{|y|^3}t}{\cos\frac{m_ek^2}{|y|^3}t}\two{x}{y}.$$  The map $D$ is called the Delaunay Hamiltonian, and $\xi_D$ is the Delaunay vector field.  Since the orbits of $\xi_D$ must preserve $|y|$, it is clear from the form of $\psi^t_D$ that $\xi_D$ is complete on $\TpS$, and every orbit is closed.

Now let $S^3_{n}$ represent $S^3$ with the north pole $(0,0,0,1)$ removed.  The Ligon-Schaaf map  $LS:N\rightarrow\TpSn$ is given by $$LS(q,p)=(A\sin\varphi+B\cos\varphi,-\nu A\cos\varphi+\nu B\sin\varphi),$$ where 
$$\begin{array}{ll}
\displaystyle\nu=\sqrt{-\frac{m_ek^2}{2H(q,p)}},&\displaystyle\varphi=\frac{1}{\nu}(q\cdot p),\\
\displaystyle A=\lb{\frac{q}{|q|}-\frac{1}{m_ek}(q\cdot p)p,\frac{1}{\nu}(q\cdot p)},&\displaystyle B=\lb{\frac{1}{\nu}|q|p,\frac{1}{m_ek}|p|^2|q|-1}.\\
\end{array}$$
This map has the following properties:
\begin{enumerate}[(1)]
\item $LS$ is a diffeomorphism between $N$ and $\TpSn$;

\item $LS^*\beta_3=\sum_{j=1}^3\lb{q_jdp_j+d(q_jp_j)}=\beta$ (recall equation \ref{eq:beta});

\item $H=D\circ LS$.
\end{enumerate}
Assertions (1) and (2) imply that $LS$ is a symplectomorphism; from that and (3), it follows that $LS_*\xi_H=\xi_D$.

\subsection{Metaplectic-c Prequantization for $\TS$}\label{subsec:mpcTS}

We begin by constructing a metaplectic-c prequantization for $(Z=T\R^4,\omega_4)$, and let it induce a metaplectic-c prequantization for $(\TS,\omega_3)$.  For $(Z,\omega_4)$, we proceed in precisely the same way as we constructed $(P,\Sigma,\gamma)$ for $(M,\omega)$ in Section \ref{subsec:mpcP}.  Choose the model vector space $V_4$, and identify the symplectic frame bundle $\Sp(Z,\omega_4)$ with $Z\times\Sp(V_4)$ using the global trivialization for the tangent bundle given by  $$\hv_j\mapsto\lat{\pp{}{x_j}}{(x,y)},\ \hw_j\mapsto\lat{\pp{}{q_j}}{(x,y)},\ \ \forall (x,y)\in Z,\ j=1,\ldots,4.$$  Let $Q_4=Z\times\Mp^c(V_4)$ with bundle projection map $Q_4\maps{\Pi_4} Z$.  Define the map $\Gamma_4:Q_4\rightarrow\Sp(T\R^4,\omega_4)$ by $$\Gamma_4(x,y,a)=(x,y,\sigma(a)),\ \ \forall (x,y)\in Z,\forall a\in\Mp^c(V_4).$$  Lastly, define the $\mathfrak{u}(1)$-valued one-form $\delta_4$ on $Q_4$ by $$\delta_4=\frac{1}{i\hbar}\Pi_4^*\beta_4+\frac{1}{2}\eta_*\vartheta_0,$$ where $\vartheta_0$ is the trivial connection on $Q_4$, and where $\beta_4$ was defined in equation \ref{eq:beta4}.  Then $(Q_4,\Gamma_4,\delta_4)$ is the unique metaplectic-c prequantization for $(Z,\omega_4)$ up to isomorphism.

In order to construct the metaplectic-c prequantization of $(\TS,\omega_3)$ induced by $(Q_4,\Gamma_4,\delta_4)$, we proceed by the process of symplectic reduction. %outlined in Section \ref{subsec:mpcquant}.  
Let $R:Z\rightarrow\R$ be given by $$R(x,y)=\frac{1}{2}|x|^2,\ \ \forall(x,y)\in Z.$$  The Hamiltonian vector field for $R$ on $Z$ is $$\xi_R=-\sum_{j=1}^4 x_j\pp{}{y_j},$$ and its flow on $Z$ is $$\psi^t_R(x,y)=(x,y-tx),\ \ \forall (x,y)\in Z.$$  It is straightforward to show that $\TS$ can be identified with the space of orbits of $\xi_R$ on the level set $A=R^{-1}(\frac{1}{2})$.  The orbit projection map $p:A\rightarrow\TS$ is $$p(x,y)=(x,y-(x\cdot y)x),\ \ \forall(x,y)\in A.$$  Let $i:A\rightarrow Z$ be the inclusion map.  Recall that $\omega_3$ is the symplectic form on $\TS$ obtained by restricting $\omega_4$.  A calculation shows that $p^*\omega_3=i^*\omega_4$.  Therefore $(\TS,\omega_3)$ is the symplectic reduction of $(Z,\omega_4)$ at the level set corresponding to $R=\frac{1}{2}$.  

Choose the model subspace $W_4$ as described in Section \ref{subsec:stdvec}, and construct the three-level structure $$(Q_{4A},\delta_{4A})\maps{\Gamma_{4A}}\Sp(TA/TA^\perp)\rightarrow A.$$  Then $\Sp(TA/TA^\perp)$ is a principal $\Sp(V_3)$ bundle over $A$, and $Q_{4A}$ is a principal $\Mp^c(V_3)$ bundle over $A$.  Let the symplectic frame bundle $\Sp(\TS,\omega_3)$ be modeled on $V_3$, so that $\Sp(\TS,\omega_3)$ is a principal $\Sp(V_3)$ bundle over $\TS$.

For each $z\in A$, $p_*|_z$ appears in the short exact sequence $$0\rightarrow T_zA^\perp\rightarrow T_zA\maps{p_*|_z}T_{p(z)}\TS\rightarrow 0.$$  Note that the projection map $p$ is the identity on $\TS\subset A$.  It follows that for all $z\in\TS$, $p_*|_z:T_zA/T_zA^\perp\rightarrow T_{z}\TS$ is a symplectic isomorphism.  Therefore $\Sp(\TS,\omega_3)$ and $\Sp(TA/TA^\perp)|_{\TS}$ are isomorphic as principal $\Sp(V_3)$ bundles over $\TS$.  Using that isomorphism, we view $\Sp(\TS,\omega_3)$ as a subbundle of $\Sp(TA/TA^\perp)$.

Let $(Q,\Gamma,\delta)$ be the result of restricting $(Q_{4A},\Gamma_{4A},\delta_{4A})$ to $\Sp(\TS,\omega_3)\subset\Sp(TA/TA^\perp)$.  Then $(Q,\Gamma,\delta)$ is a metaplectic-c prequantization for $(\TS,\omega_3)$, and it is unique up to isomorphism.  We will also use the notation $(Q,\Gamma,\delta)$ to denote the metaplectic-c prequantizations for $\TpS$ and $\TpSn$ obtained by restriction.

\subsection{Lifting the Ligon-Schaaf Map}

Recall that the Ligon-Schaaf map $LS:N\rightarrow\TpSn$ is a symplectomorphism.  Define the map $\widetilde{LS}:\Sp(N,\omega)\rightarrow\Sp(\TpSn,\omega_3)$ by $$\widetilde{LS}(b)=LS_*\circ b,\ \ \forall b\in\Sp(N,\omega),$$ and observe that this is an isomorphism of principal $\Sp(V_3)$ bundles.  

In Section \ref{subsec:mpcP}, we used the global trivialization $$\hv_j\mapsto\lat{\pp{}{q_j}}{m},\ \ \hw_j\mapsto\lat{\pp{}{p_j}}{m},\ \ \forall m\in N,\ j=1,2,3,$$ to identify $\Sp(N,\omega)$ with $N\times\Sp(V_3)$.  From the map $\widetilde{LS}$, we see that $\Sp(\TpSn,\omega_3)$ can be identified with $\TpSn\times\Sp(V_3)$ with respect to the global trivialization $$\hv_j\mapsto\lat{LS_*\pp{}{q_j}}{LS(m)},\ \ \hw_j\mapsto\lat{LS_*\pp{}{q_j}}{LS(m)},\ \ \forall LS(m)\in\TpSn,\ j=1,2,3.$$  In terms of these trivializations for $\Sp(N,\omega)$ and $\Sp(\TpSn,\omega_3)$, $\widetilde{LS}$ is given simply by 
\begin{equation}\label{eq:lstilde}
\widetilde{LS}(m,g)=(LS(m),g),\ \ \forall m\in N,\forall g\in\Sp(V_3).
\end{equation}

As we have seen before, once we have a global trivialization for the symplectic frame bundle, it is straightforward to construct a metaplectic-c prequantization.  Let $Q'=\TpSn\times\Mp^c(V_3)$ with bundle projection map $Q'\maps{\Pi'}\TpSn$.  Define the map $\Gamma':\TpSn\rightarrow\Sp(\TpSn,\omega_3)$ by $$\Gamma'(z,a)=(z,\sigma(a)),\ \ \forall z\in\TpSn,\forall a\in\Mp^c(V_3).$$  Let  $$\delta'=\frac{1}{i\hbar}\Pi'^*\beta_3+\frac{1}{2}\eta_*\vartheta_0,$$ where $\beta_3$ is the restriction to $\TpSn$ of the one-form defined in equation \ref{eq:beta4}, and where $\vartheta_0$ is the trivial connection on $Q'$.  Then $(Q',\Gamma',\delta')$ is a metaplectic-c prequantization for $(\TpSn,\omega_3)$.  Since $\TpSn$ is simply connected, $(Q',\Gamma',\delta')$ must be isomorphic to $(Q,\Gamma,\delta)$.  

Recall that $P=M\times\Mp^c(V_3)$.  Define $\widehat{LS}:P\rightarrow Q'$ by $$\widehat{LS}(m,a)=(LS(m),a).$$  This is clearly an isomorphism of principal $\Mp^c(V_3)$ bundles.    Since $LS^*\beta_3=\beta$, we have $LS^*\delta'=\gamma$.  Lastly, it follows from equation \ref{eq:lstilde} that $\Gamma'\circ\widehat{LS}=\widetilde{LS}\circ\Sigma$.  Therefore $\widehat{LS}$ is an isomorphism of metaplectic-c prequantizations.

All of the preceding observations combine to yield the following commutative diagram.
$$\xymatrix{(P,\gamma) \ar[d]^{\Sigma} \ar[r]^(0.45){\widehat{LS}} & (Q',\delta') \ar[d]^{\Gamma'}\\
\Sp(N,\omega) \ar[d] \ar[r]^(0.45){\widetilde{LS}} & \Sp(\TpSn,\omega_3) \ar[d]\\
(N,\omega) \ar[r]^(0.45){LS}& (\TpSn,\omega_3)
}$$
Each of the maps $LS$, $\widetilde{LS}$, and $\widehat{LS}$ is an isomorphism.  From these isomorphisms and the fact that $D\circ LS=H$, it follows that $E<0$ is a quantized energy level of $(N,\omega,H)$ if and only if it is a quantized energy level of $(\TpSn,\omega_3,D)$.  

In fact, we claim that $E<0$ is a quantized energy level of $(N,\omega,H)$ if and only if it is a quantized energy level of $(\TpS,\omega_3,D)$.  When we replace the north pole in $S^3$, we acquire more closed orbits:  namely, those with $x$-components that pass through $(1,0,0,0)$.  However, due to the rotational symmetry of the system $(\TpS,\omega_3,D)$, an orbit that passes through $(1,0,0,0)$ can always be transformed into one that does not, without altering the holonomy condition.  Therefore the quantized energy levels of $(\TpSn,\omega_3,D)$ and $(\TpS,\omega_3,D)$ are identical.

\subsection{Rescaling the Delaunay Hamiltonian}\label{subsec:rescale}

So far, we have shown that the negative quantized energy levels of the hydrogen atom are the same as the quantized energy levels of the Delaunay Hamiltonian on $\TpS$.  In this section, we will make one final transformation that relates these energies to the quantized energy levels of a free particle on $\TS$.

Let $K:TS^3\rightarrow\R$ be given by $$K(x,y)=\frac{1}{2}|y|^2,\ \ \forall(x,y)\in\TS.$$  In the context of classical mechanics, this energy function describes a free particle on $\TS$.  It is shown in \cite{cb} that the corresponding Hamiltonian vector field on $\TS$ is $$\xi_K=\sum_{j=1}^4\lb{y_j\pp{}{x_j}-|y|^2x_j\pp{}{y_j}},$$ and the flow of this vector field is 
\begin{equation}\label{eq:flowK}
\psi^t_K(x,y)=\four{\cos|y|t}{\frac{1}{|y|}\sin|y|t}{-|y|\sin|y|t}{\cos|y|t}\two{x}{y},\ \ \forall(x,y)\in\TS\ \mbox{such that}\ |y|>0.
\end{equation}
If $|y|=0$, then the particle is stationary, and the flow is simply $\psi^t_K(x,0)=(x,0)$.

On the submanifold $\TpS$, the range of $K$ is $\R_{>0}$, while the range of the Delaunay Hamiltonian $D$ is $\R_{<0}$.  Let $f:\R_{>0}\rightarrow\R_{<0}$ be given by 
\begin{equation}\label{eq:diffeo}
f(z)=-\frac{m_ek^2}{4z},\ \ \forall z\in\R_{>0}.
\end{equation}
Note that $f$ is a diffeomorphism, and $D=f\circ K$ on $\TpS$.  Using the dynamical invariance property of quantized energy levels that was stated in Section \ref{subsec:mpcquant}, we see that $E$ is a quantized energy level of $(\TpS,\omega_3,D)$ if and only if $f^{-1}(E)$ is a quantized energy level of $(\TpS,\omega_3,K)$.  Thus the quantized energy levels of $(N,\omega,H)$ are exactly the images of the positive quantized energy levels of $(\TS,\omega_3,K)$ under the diffeomorphism $f$.

\section{Quantization of a Free Particle on $\TS$}\label{sec:freeTS}

\subsection{Orbits of $\xi_K$ and Local Coordinates for $\TS$}

Let $\E>0$ be arbitrary, and let $\sS=K^{-1}(\E)\subset\TS$.  From the flow of $\xi_K$ in equation \ref{eq:flowK}, it is apparent that all orbits of $\xi_K$ in $\sS$ are closed, with period $\frac{2\pi}{\sqrt{2\E}}$.  Let $z_0=(x_0,y_0)\in\sS$ be an arbitrary initial point, and let $\sC\subset\sS$ be the the orbit of $\xi_K$ through $z_0$:   $$\sC=\cb{\psi^t_K(z_0):t\in\R}\subset\TS,$$ 

We can use the rotational symmetry of the system $(\TS,\omega_3,K)$ to make some simplifying assumptions about $\sC$.  If we view $x_0$ and $y_0$ as two perpendicular vectors in $\R^4$, then there is some rotation about the origin in $\R^4$ that carries them both to the $x_3x_4$-plane.  By performing this rotation in $x$-space and $y$-space, we can assume without loss of generality that $x_0$ and $y_0$ take the form $x_0=(0,0,x_{30},x_{40})$ and $y_0=(0,0,y_{30},y_{40})$, where $|x|^2=x_{30}^2+x_{40}^2=1$, $x\cdot y=x_{30}y_{30}+x_{40}y_{40}=0$, and $|y|^2=y_{30}^2+y_{40}^2=2k$.  The orbit $\sC$ then lies in $x_3x_4y_3y_4$-space.

Thus far, we have treated $\TS$ as a submanifold of $Z=T\R^4$, using the coordinates $(x,y)$ on $Z$ to describe points in $\TS$.  Now we make a local change of coordinates on $Z$ that will yield symplectic coordinates for $\TS$ on a neighborhood that contains $\sC$.  Specifically, we introduce 4-dimensional spherical coordinates and their conjugate momenta.  Let $U\subset Z$ be the open set $$U=\cb{(x,y)\in T\R^4:x_3^2+x_4^2>0}.$$  On $U$, let the new spatial coordinates be $(a,b,c,r)$, where
\begin{align*}
a=&\arctan\frac{\sqrt{x_2^2+x_3^2+x_4^2}}{x_1},\\
b=&\arctan\frac{\sqrt{x_3^2+x_4^2}}{x_2},\\
c=&\arctan\frac{x_4}{x_3},\\
r=&\sqrt{x_1^2+x_2^2+x_3^2+x_4^2}.
\end{align*}
The angles $a$ and $b$ are defined modulo $\pi$, and the angle $c$ is defined modulo $2\pi$.  Let $\rho_1=\sqrt{x_2^2+x_3^2+x_4^2}$ and $\rho_2=\sqrt{x_3^2+x_4^2}$.  The conjugate momenta corresponding to the spherical coordinates are  
\begin{align*}
p_a=&\frac{x_1}{\rho_1}(x_2y_2+x_3y_3+x_4y_4)-\rho_1y_1,\\
p_b=&\frac{x_2}{\rho_2}(x_3y_3+x_4y_4)-\rho_2y_2,\\
p_c=&x_3y_4-x_4y_3,\\
p_r=&\frac{1}{r}(x_1y_1+x_2y_2+x_3y_3+x_4y_4).
\end{align*}
Later, we will need the inverse transformations, which are
\begin{equation}\label{eq:convert1}
\begin{array}{l}
x_1=r\cos a,\\
x_2=r\sin a\cos b,\\
x_3=r\sin a\sin b\cos c,\\
x_4=r\sin a\sin b\sin c,\\
\end{array}
\end{equation}

\begin{equation}\label{eq:convert2}
\begin{array}{l}
\displaystyle y_1=p_r\cos a-\frac{p_a}{r}\sin a,\\
\displaystyle y_2=p_r\sin a\cos b+\frac{p_a}{r}\cos a\cos b-\frac{p_b}{r\sin a}\sin b\\
\displaystyle y_3=p_r\sin a\sin b\cos c+\frac{p_a}{r}\cos a\sin b\cos c+\frac{p_b}{r\sin a}\cos b\cos c-\frac{p_c}{r\sin a\sin b}\sin c,\\
\displaystyle y_4=p_r\sin a\sin b\sin c+\frac{p_a}{r}\cos a\sin b\cos c+\frac{p_b}{r\sin a}\cos b\sin c+\frac{p_c}{r\sin a\sin b}\cos c.
\end{array}
\end{equation}
For convenience, we let $a_j$, $j=1,\ldots,4$, range over $a,b,c,r$.  

On $U$, one can verify that 
\begin{equation}\label{eq:beta4s}
\beta_4=-\sum_{j=1}^4y_jdx_j=-\sum_{j=1}^4p_{a_j}da_j,
\end{equation}
and so $$\omega_4=\sum_{j=1}^4dx_j\wedge dy_j=\sum_{j=1}^4 da_j\wedge dp_j.$$  The submanifold $\TS$ is characterized by the constant values $r=1$ and $p_r=0$, which implies that the restrictions to $\TS\cap U$ of $\beta_4$ and $\omega_4$ are $$\beta_3=-\sum_{j=1}^3 p_{a_j}da_j,\ \ \omega_3=\sum_{j=1}^3da_j\wedge dp_{a_j}.$$  Thus $(a,b,c,p_a,p_b,p_c)$ are symplectic coordinates for $\TS\cap U$.  On this neighborhood, the map $K$ takes the form  $$K=\frac{1}{2}\lb{p_a^2+\frac{p_b^2}{\sin^2a}+\frac{p_c^2}{\sin^2a\sin^2b}}.$$

Several times now, once we had a set of symplectic coordinates such as $(a,b,c,r,p_a,p_b,p_c,p_r)$, we used the trivialization of the symplectic frame bundle given by the coordinate vector fields to construct a metaplectic-c prequantization.  We could apply this same procedure to $U$; however, since $U$ is not simply connected, it is not necessarily the case that the metaplectic-c prequantization so constructed would be isomorphic to the result of restricting $(Q_4,\Gamma_4,\delta_4)$ to $U$.  Instead, we must show how the local change of variables from Cartesian to spherical coordinates can be lifted from $Z$ to $Q_4$.  This is the subject of the next section.

\subsection{Change of variables over $\sC$}

Recall that the integral curve of $\xi_K$ through the initial point $z_0=(x_0,y_0)$ is $$\psi^t_K(x_0,y_0)=\four{\cos\sqrt{2\E}t}{\frac{1}{\sqrt{2\E}}\sin\sqrt{2E}t}{-\sqrt{2\E}\sin\sqrt{2\E}t}{\cos\sqrt{2\E}t}\two{x_0}{y_0},$$ and its image $\sC$ is a closed curve lying in $x_3x_4y_3y_4$-space.  Since $\sC\subset U$, the points in $\sC$ can be rewritten in spherical coordinates.  Upon converting, we find that any point in $\sC$ satisfies $p_a=p_b=0$ and $a=b=\frac{\pi}{2}$.  Further, $p_c$ is a constant value over $\sC$ satisfying $p_c^2=2\E$.  Since $\sC\subset\TS$, we also have $r=1$ and $p_r=0$.  Therefore, in spherical coordinates, the orbit takes the form $$\sC=\cb{\lb{\frac{\pi}{2},\frac{\pi}{2},c,1,0,0,p_c,0}:c\in\R/2\pi\Z}.$$  Let $z(c)$ represent the point $\lb{\frac{\pi}{2},\frac{\pi}{2},c,1,0,0,p_c}\in\sC$.

The change of coordinates from Cartesian to spherical must now be lifted to the symplectic frame bundle and the metaplectic-c prequantization for $(Z,\omega_4)$.  The change of coordinates on $\Sp(Z,\omega_4)$ will take place over $U$, and that on $Q_4$ will take place over $\sC$.  This procedure was first demonstrated in \cite{v1}, and we follow the same steps now.

On the neighborhood $U$, we have two different coordinate maps:  the Cartesian map $\Phi_c:U\rightarrow\R^8$, and the spherical map $\Phi_s:U\rightarrow \R^4\times(\R/\pi\Z)^2\times\R/2\pi\Z\times\R$.  The change of variables on $U$ is simply the transition map $F=\Phi_s\circ\Phi_c^{-1}$.  Let $\Phi_c(U)=U_c$ and $\Phi_s(U)=U_s$.  Then each of the following maps is a diffeomorphism.
\begin{center}
$\xymatrix{
& U \ar[ld]_{\Phi_c} \ar[rd]^{\Phi_s} & \\
U_c \ar[rr]^F & & U_s
}$
\end{center}
These observations are straightforward on $U$, but they will motivate the constructions on the symplectic frame bundle and the metaplectic-c prequantization.

Let $b_c$ be the section of $\Sp(Z,\omega_4)$ over $U$ given by $$b_c(z):V_4\rightarrow T_{z}Z\ \ \mbox{such that}\ \ \hv_j\mapsto\lat{\pp{}{x_j}}{z},\ \hw_j\mapsto\lat{\pp{}{y_j}}{z},\ \ \forall z\in U,\ j=1,\ldots,4.$$  That is, $b_c$ is the section that defines the trivialization of $\Sp(Z,\omega_4)|_U$ with respect to Cartesian coordinates.  Let the map $\widetilde{\Phi}_c:\Sp(Z,\omega)|_{U}\rightarrow U_c\times\Sp(V_4)$ be given by $$\widetilde{\Phi}_c(b_c(z)\cdot g)=(\Phi_c(z),g),\ \ \forall z\in U,\ \forall g\in\Sp(V_4).$$  Similarly, let $b_s$ be the section of $\Sp(Z,\omega_4)|_U$ given by 
$$b_s(z):V_4\rightarrow T_zZ\ \ \mbox{such that}\ \ \hv_j\mapsto\lat{\pp{}{a_j}}{z},\ \hw_j\mapsto\lat{\pp{}{p_{a_j}}}{z},\ \ \forall z\in U,\ j=1,\ldots,4,$$ and define the map $\widetilde{\Phi}_s:\Sp(Z,\omega)|_U\rightarrow U_s\times\Sp(V_4)$ by
$$\widetilde{\Phi}_s(b_s(z)\cdot g)=(\Phi_s(z),g),\ \ \forall z\in U,\ \forall g\in\Sp(V_4).$$  

To perform the change of coordinates on the level of the symplectic frame bundle, we must lift $F$ to a map $\widetilde{F}:U_c\times\Sp(V_4)\rightarrow U_s\times\Sp(V_4)$ in such a way that the following diagram commutes.  
\begin{center}
$\xymatrix{
& \Sp(Z,\omega)|_{U} \ar[ld]_{\widetilde{\Phi}_c} \ar[rd]^{\widetilde{\Phi}_s} & \\
U_c\times\Sp(V_4) \ar[rr]^{\widetilde{F}} & & U_s\times\Sp(V_4)
}$
\end{center}

In equations \ref{eq:convert1} and \ref{eq:convert2}, we gave explicit formulas for $x_j$ and $y_j$ in terms of $a_k$ and $p_{a_k}$.  At each $z\in U$, let $G(z)$ be the $8\times 8$ matrix consisting of the partial derivatives of the Cartesian coordinates with respect to the spherical ones:  
$$G(z)=\four{\lat{\pp{x_k}{a_j}}{z}}{\lat{\pp{y_k}{a_j}}{z}}{\lat{\pp{x_k}{p_{a_j}}}{z}}{\lat{\pp{y_k}{p_{a_j}}}{z}}_{1\leq j,k\leq 4}.$$  Then $$G(z)\two{\lat{\pp{}{x_k}}{z}}{\lat{\pp{}{y_k}}{z}}_{1\leq k\leq 4}=\two{\lat{\pp{}{a_j}}{z}}{\lat{\pp{}{p_{a_j}}}{z}}_{1\leq j\leq 4},\ \ \forall z\in U.$$  Since the Cartesian and spherical coordinate vectors are both symplectic bases for $T_zZ$, $G(z)$ is a symplectic matrix and can be treated as an element of $\Sp(V_4)$.  We claim that the desired map $\widetilde{F}$ is given by $$\widetilde{F}(\Phi_c(z),g)=(\Phi_s(z),G(z)g),\ \ \forall z\in U,\ \forall g\in\Sp(V_4).$$ 

To prove this claim, it suffices to show that $$\widetilde{\Phi}_s(b_c(z))=\widetilde{F}\circ\widetilde{\Phi}_c(b_c(z));$$ the images under $\widetilde{F}$ of all other elements of the fiber $\Sp(Z,\omega)_z$ are determined by the group action.  On the right-hand side, we have $$\widetilde{F}\circ\widetilde{\Phi}_c(b_c(z))=\widetilde{F}(\Phi_c(z),I)=(\Phi_s(z),G(z)).$$  Using the definition of $b_c$, we find that $$b_c(z):G(z)\two{\hv_j}{\hw_j}_{1\leq j\leq 4}\mapsto G(z)\two{\lat{\pp{}{x_j}}{z}}{\lat{\pp{}{y_j}}{z}}_{1\leq j\leq 4}=\two{\lat{\pp{}{a_j}}{z}}{\lat{\pp{}{p_{a_j}}}{z}}_{1\leq j\leq 4}.$$  Therefore, with respect to the spherical trivialization, $$\widetilde{\Phi}_s(b_c(z))=(\Phi_s(z),G(z))=\widetilde{F}\circ\widetilde{\Phi}_c(b_c(z)),$$ as needed.  Hence the map $\widetilde{F}$ effects the change of coordinates from Cartesian to spherical on $\Sp(Z,\omega_4)|_{U}$.

The final lift to the metaplectic-c prequantization will take place over $\sC$.  Let $\sC_c$ and $\sC_s$ be the images of $\sC$ under $\Phi_c$ and $\Phi_s$, respectively.  Recall that points in $\sC$ are denoted by $z(c)$ with $c\in\R/2\pi\Z$.  We write $G(c)$ for $G(z(c))$.  The components of $G(c)$ are single-valued with respect to $c$, so $G(c)$ is a closed path through $\Sp(V_4)$.

It is clear from the construction of $Q_4$ how to define a local trivialization with respect to the Cartesian coordinates.  Let $\widehat{\Phi}_c:Q_4|_{\sC}\rightarrow\sC_c\times\Mp^c(V_4)$ be given by  $\widehat{\Phi}_c(z(c),a)=(\Phi_c(z(c)),a)$ for all $z(c)\in\sC$ and $a\in\Mp^c(V_4)$.  Then $\widehat{\Phi}_c$ and $\widetilde{\Phi}_c$ make the following diagram commute.
\begin{center}
$\xymatrix{
Q_4|_{\sC} \ar[r]^{\Gamma_4} \ar[d]^{\widehat{\Phi}_c}  & \Sp(Z,\omega_4)|_{\sC} \ar[d]^{\widetilde{\Phi}_c} \\
\sC_{c}\times\Mp^c(V_4) \ar[r]^{\sigma} & \sC_c\times\Sp(V_4)
}$
\end{center}
We require a local trivialization of $Q_4$ with respect to the spherical coordinates that has the analogous relationship to $\widetilde{\Phi}_s$.  To find the appropriate map $\widehat{\Phi}_s:Q_4|_{\sC}\rightarrow\sC_s\times\Mp^c(V_4)$, we will construct the map $\widehat{F}$ shown below, then define $\widehat{\Phi}_s=\widehat{F}\circ\widehat{\Phi}_c$.
\begin{center}
$\xymatrix{
& Q_4|_{\sC} \ar[dl]_{\widehat{\Phi}_c} \ar[rd] \ar[rrrd]^{\Sigma} & \\
\sC_c\times\Mp^c(V_4) \ar[rr]^{\widehat{F}} \ar[rrrd]^{\sigma} & & \sC_s\times\Mp^c(V_4) \ar[rrrd]^(.25){\sigma} & & \Sp(Z,\omega_4)|_{\sC} \ar[ld]_(.4){\widetilde{\Phi}_c} \ar[rd]^(.4){\widetilde{\Phi}_s} & \\
& & & \sC_c\times\Sp(V_4) \ar[rr]^{\widetilde{F}} & & \sC_s\times\Sp(V_4)
}$
\end{center} 
An examination of the condition $\sigma\circ\widehat{F}\circ\widehat{\Phi}_c=\widetilde{F}\circ\sigma\circ\widehat{\Phi}_c$ shows that if $\widehat{F}$ exists, then it has the form $$\widehat{F}(\Phi_c(z(c)),a)=(\Phi_s(z(c)),\widehat{G}(c)a),\ \ \forall z(c)\in\sC,\forall a\in\Mp^c(V_4),$$ where $\widehat{G}(c)\in\Mp^c(V_4)$ and $\sigma(\widehat{G}(c))=G(c)$ for all $c$.  That is, $\widehat{G}(c)$ must be a lift of the path $G(c)$ to $\Mp^c(V_4)$.  Further, for $\widehat{F}$ to be a well-defined function, $\widehat{G}(c)$ must be single-valued with respect to $c$.

Since $\Mp^c(V_4)$ is a circle extension of $\Sp(V_4)$, there are many possible lifts of $G(c)$ to $\Mp^c(V_4)$.  However, it would be ideal if $\widehat{G}(c)$ were a path not only through $\Mp^c(V_4)$ but through the subgroup $\Mp(V_4)$.  Recall that the one-form $\delta_4$ on $Q_4$ is $$\delta_4=\frac{1}{i\hbar}\Pi^*\beta_4+\frac{1}{2}\eta_*\vartheta_0,$$ where $\vartheta_0$ is the trivial connection on $Z\times\Mp^c(V_4)$.  Since $\ker\eta=\Mp(V_4)$, a map $\widehat{F}$ that involves only $\widehat{G}(c)\in\Mp(V)$ would not alter the form of the $\frac{1}{2}\eta_*\vartheta_0$ term in $\delta_4$.  Therefore we will lift $G(c)$ to $\widehat{G}(c)\in\Mp(V_4)$, and determine whether the lifted path satisfies $\widehat{G}(c+2\pi)=\widehat{G}(c)$.  

If $\widehat{G}(c)\in\Mp(V_4)$, then its parameters take the form $(G(c),\mu(c))$, where $\mu(c)^2\mbox{Det}_\C C_{G(c)}=1$.  To determine $\mu(c)$, we must calculate $G(c)$.  Using the expressions in equations \ref{eq:convert1} and \ref{eq:convert2}, we compute the partial derivatives that form the matrix $G$, and evaluate at the point $z(c)=\lb{\frac{\pi}{2},\frac{\pi}{2},c,1,0,0,p_c,0}\in\sC_s$.  Using the abbreviations $S(c)=\sin c$ and $C(c)=\cos c$, the result is
$$G(c)=\lb{\begin{array}{cccccccc}
-1&0&0&0&0&0&0&0\\
0&-1&0&0&0&0&0&0\\
0&0&-S(c)&C(c)&0&0&-p_cC(c)&-p_cS(c)\\
0&0&C(c)&S(c)&0&0&p_cS(c)&-p_cC(c)\\
0&0&0&0&-1&0&0&0\\
0&0&0&0&0&-1&0&0\\
0&0&0&0&0&0&-S(c)&C(c)\\
0&0&0&0&0&0&C(c)&S(c)
\end{array}}.$$

Next, we evaluate $C_{G(c)}=\frac{1}{2}(G(c)-JG(c)J)$, using the matrix form for $J$ noted in Section \ref{subsec:stdvec}, and convert it to a $4\times 4$ complex matrix.  We find that  $$C_{G(c)}=\lb{\begin{array}{cccc}
-1&0&0&0\\
0&-1&0&0\\
0&0&-S(c)-\frac{i}{2}p_cC(c)&C(c)-\frac{i}{2}p_cS(c)\\
0&0&C(c)+\frac{i}{2}p_cS(c)&S(c)-\frac{i}{2}p_cC(c)
\end{array}},$$ which has complex determinant $$\Det_\C C_{G(c)}=-1-\frac{1}{4}p_c^2.$$  This value is real and constant over $\sC$.  Thus we can define $\widehat{G}(c)$ to be the element of $\Mp(V_4)\subset\Mp^c(V_4)$ with parameters $\lb{G(c),i\lb{1+\frac{1}{4}p_c^2}^{-1/2}}$, for all $c\in\R/2\pi\Z$, whereupon $\widehat{G}(c)$ is the desired closed path through $\Mp(V_4)$.

Having determined $\widehat{G}(c)$, we now define the map $\widehat{F}:\sC_c\times\Mp^c(V_4)\rightarrow\sC_s\times\Mp^c(V_4)$ by $$\widehat{F}(\Phi_c(z(c)),a)=(\Phi_s(z(c)),\widehat{G}(c)a)\ \ \forall z(c)\in\sC,\forall a\in\Mp^c(V_4).$$  The local trivialization of $Q_4|_{\sC}$ that is compatible with spherical coordinates comes about by setting $\widehat{\Phi}_s=\widehat{F}\circ\widehat{\Phi}_c$.  The one-form $\delta_4$ on $Q_4$ induces a one-form $\delta_{4s}$ on $\sC_s\times\Mp^c(V_4)$ that takes the form $$\delta_{4s}=\frac{1}{i\hbar}\Pi^*\beta_4+\frac{1}{2}\eta_*\vartheta_0,$$ where $\beta_4$ is written in spherical coordinates as in equation \ref{eq:beta4s}, and where $\vartheta_0$ is now the trivial connection on $\sC_s\times\Mp^c(V_4)$.

\subsection{Restrictions to $\sC\subset\TS$}

From now on, we no longer write the maps $\Phi_s$, $\widetilde{\Phi}_s$ and $\widehat{\Phi}_s$ explicitly, but write elements of $Z$, $\Sp(Z,\omega_4)|_U$ and $Q_4|_{\sC}$ with respect to spherical coordinates and the local spherical trivializations.  Further, we treat $\TS$ as a six-dimensional manifold with local coordinates $(a,b,c,p_a,p_b,p_c)$ on $\TS\cap U$, and with symplectic form $\omega_3=\sum_{j=1}^3 da_j\wedge dp_{a_j}$.

Recall from Section \ref{subsec:mpcTS} that we defined the map $R(x,y)=\frac{1}{2}|x|^2$ on $Z$, and set $A=R^{-1}(\frac{1}{2})$.  We argued that the symplectic frame bundle $\Sp(\TS,\omega_3)$ is isomorphic to $\Sp(TA/TA^\perp)|_{\TS}$, and that the metaplectic-c prequantization $(Q,\Gamma,\delta)$ for $(\TS,\omega_3)$ is obtained by restricting $(Q_{4A},\Gamma_{4A},\delta_{4A})$ to $\Sp(TA/TA^\perp)|_{\TS}$.  

In spherical coordinates, we have $$R(z)=\frac{1}{2}r^2,\ \ \forall z\in Z\cap U,$$ which has Hamiltonian vector field $$\xi_R=r\pp{}{p_r}.$$  It is immediate that at each point $z\in A\cap U$, $$T_zA=\mbox{span}\lat{\cb{\pp{}{a},\pp{}{b},\pp{}{c},\pp{}{p_a},\pp{}{p_b},\pp{}{p_c},\pp{}{p_r}}}{z}\ \ \mbox{and}\ \ T_zA^\perp=\mbox{span}\lat{\cb{\pp{}{p_r}}}{z},$$ which implies that $$T_zA/T_zA^\perp=\mbox{span}\lat{\cb{\ls{\pp{}{a}},\ls{\pp{}{b}},\ls{\pp{}{c}},\ls{\pp{}{p_a}},\ls{\pp{}{p_b}},\ls{\pp{}{p_c}}}}{z}.$$  We see that over $A\cap U$, the local trivialization of $\Sp(Z,\omega_4)|_U$ with respect to spherical coordinates induces a local trivialization $$\Sp(Z,\omega_4;A)|_{A\cap U}=A\cap U\times\Sp(V_4;W_4).$$  The group homomorphism $\nu:\Sp(V_4;W_4)\rightarrow\Sp(V_3)$ then induces a local trivialization $$\Sp(TA/TA^\perp)|_{A\cap U}=A\cap U\times\Sp(V_3).$$  Restricting further to $\TS\cap U\subset A\cap U$, we find that $$\Sp(\TS,\omega_3)|_{\TS\cap U}=\TS\cap U\times\Sp(V_3),$$ with the local trivialization given by $$\hv_j\mapsto\lat{\pp{}{a_j}}{z},\ \ \hw_j\mapsto\lat{\pp{}{p_{a_j}}}{z},\ \ \forall z\in\TS\cap U,\ j=1,2,3.$$

On the level of metaplectic-c bundles, we can write $$Q_4^A|_{\sC}=\sC\times\Mp^c(V_4;W_4)$$ by restricting $Q_4|_{\sC}$ to $\Sp(Z,\omega;A)|_{\sC}$.  Then the group homomorphism $\hat{\nu}:\Mp^c(V_4;W_4)\rightarrow\Mp^c(V_3)$ induces the local trivialization $$Q_{4A}|_{\sC}=Q|_{\sC}=\sC\times\Mp^c(V_3).$$  The prequantization one-form on $Q$ is $$\delta=\frac{1}{i\hbar}\Pi^*\beta_3+\frac{1}{2}\eta_*\vartheta_0,$$ where $\vartheta_0$ is the trivial connection on $\sC\times\Mp^c(V_3)$, and where $$\beta_3=-\sum_{j=1}^3 p_jda_j=-p_cdc$$ over $\sC$.

\subsection{Quantized Energy Levels for $(\TS,\omega_3,K)$}

In the local coordinates $(a,b,c,p_a,p_b,p_c)$, the map $K:\TS\rightarrow\R$ is given by  $$K=\frac{1}{2}\lb{p_a^2+\frac{p_b^2}{\sin^2a}+\frac{p_c^2}{\sin^2a\sin^2b}},$$ and its Hamiltonian vector field is  $$\xi_K=p_a\pp{}{a}+\frac{p_b}{\sin^2a}\pp{}{b}+\frac{p_c}{\sin^2a\sin^2b}\pp{}{c}+\lb{\frac{p_b^2}{\sin^3 a}\cos a+\frac{p_c^2\cos a}{\sin^3 a\sin^2b}}\pp{}{p_a}+\frac{p_c^2\cos b}{\sin^2a\sin^3b}\pp{}{p_b}.$$  Let $\xi_K$ have flow $\psi_K^t$ on $\TS$.  

Recall that we fixed $\E>0$ and defined $\sS=K^{-1}(\E)$, and that $z_0$ is a point in $\sS$.  The closed curve $\sC$ is the orbit of $\xi_K$ through $z_0$.  We will show that the orbit of $\txi_K$ through $(z_0,I)\in\Sp(T\sS/T\sS^\perp)$ is closed, and then we will determine the values of $\E$ for which the orbit of $\hxi_K$ through $(z_0,I)\in Q_\sS$ is also closed.

The lift of $\psi^t_K$ to $\widetilde{\psi}^t_K$ on $\Sp(\TS,\omega_3)$ is given by $$\widetilde{\psi}^t_K(z,I)=(z,\psi^t_{K*}|_z),\ \ \forall z\in\TS\cap U,$$ which implies that the lifted vector field $\txi_K$ is $$\txi_{K}(z,I)=\lb{\xi_K(z),\lat{\dd{}{t}}{t=0}\psi^t_{K*}|_z},\ \ \forall z\in\TS\cap U.$$  As an $8\times 8$ matrix, $\lat{\dd{}{t}}{t=0}\psi^t_{K*}|_z$ can be interpreted as an element of the Lie algebra $\mathfrak{sp}(V_3)$, and its components are  $$\lb{\lat{\dd{}{t}}{t=0}\psi^t_*|_{z}}_{jk}=\lat{\pp{\lb{\xi_K}_j}{X_k}}{X=z},$$ where $X_k$ ranges over $(a,b,c,p_a,p_b,p_c)$.  In particular, when we evaluate this matrix of partial derivatives at $z(c)\in\sC$, we obtain $$\lat{\dd{}{t}}{t=0}\psi^t_{K*}|_{z(c)}=
\lb{\begin{array}{cccccc}
0&0&0&1&0&0\\
0&0&0&0&1&0\\
0&0&0&0&0&1\\
-p_c^2&0&0&0&0&0\\
0&-p_c^2&0&0&0&0\\
0&0&0&0&0&0\end{array}}.$$  
This value is constant over $\sC$; we denote it by $\kappa$.  Thus $$\txi_K(z(c),I)=(\xi_K(z(c)),\kappa),\ \ \forall z(c)\in\sC,$$ which implies that the flow of $\txi_K$ through $(z_0,I)\in\Sp(\TS,\omega_3)$ is $$\widetilde{\psi}_K^t(z_0,I)=(\psi_K^t(z_0),\exp(t\kappa)).$$  Let $\lambda=\sqrt{p_c^2}=\sqrt{2\E}$.  A calculation shows that $$\exp(t\kappa)=
\lb{\begin{array}{cccccc}
\cos\lambda t&0&0&\frac{1}{\lambda}\sin\lambda t&0&0\\
0&\cos\lambda t&0&0&\frac{1}{\lambda}\sin\lambda t&0\\
0&0&1&0&0&0\\
-\lambda\sin\lambda t&0&0&\cos\lambda t&0&0\\
0&-\lambda\sin\lambda t&0&0&\cos\lambda t&0\\
0&0&0&0&0&1\end{array}}.$$  
Thus the orbit through $(x_0,I)$ is closed, with period $\frac{2\pi}{\lambda}$.  

Over $\sC$, $\xi_K$ reduces to $$\xi_K=p_c\pp{}{c}.$$  Therefore, for all $z(c)\in\sC$, we have  $$T_{z(c)}\sS^\perp=\mbox{span}\lat{\cb{\pp{}{c}}}{z(c)},\ \ T_{z(c)}\sS=\mbox{span}\lat{\cb{\pp{}{a},\pp{}{b},\pp{}{c},\pp{}{p_a},\pp{}{p_b}}}{z(c)},$$ 
and  $$T_{z(c)}\sS/T_{z(c)}\sS^\perp=\mbox{span}\lat{\cb{\ls{\pp{}{a}},\ls{\pp{}{b}},\ls{\pp{}{p_a}},\ls{\pp{}{p_b}}}}{z(c)}.$$  We identify $\Sp(T\sS/T\sS^\perp)|_{\sC}$ with $\sC\times\Sp(V_2)$ in the obvious way.  

Upon descending to $\Sp(T\sS/T\sS^\perp)|_{\sC}$, the induced flow takes the form $$\widetilde{\psi}_K^t(z_0,I)=(\psi^t_K(z_0),\nu(\exp(t\kappa))),$$ where $\nu$ is the group homomorphism $\Sp(V_3;W_3)\maps{\nu}\Sp(V_2)$.  We calculate that $$\nu(\exp(t\kappa))=\lb{\begin{array}{cccc}\cos\lambda t&0&\frac{1}{\lambda}\sin\lambda t&0\\0&\cos\lambda t&0&\frac{1}{\lambda}\sin\lambda t\\-\lambda\sin\lambda t&0&\cos\lambda t&0\\0&-\lambda\sin\lambda t&0&\cos\lambda t\end{array}}.$$  The corresponding vector field on $\Sp(T\sS/T\sS^\perp)|_{\sC}$ is $$\txi_K(z(c),I)=(\xi_K(z(c)),\overline{\kappa}),\ \ \forall z(c)\in\sC,$$ where $$\overline{\kappa}=\lat{\dd{}{t}}{t=0}\nu(\exp(t\kappa))=\lb{\begin{array}{cccc}0&0&1&0\\0&0&0&1\\-\lambda^2&0&0&0\\0&-\lambda^2&0&0\end{array}}\in\mathfrak{sp}(V_2).$$  

The local trivializations lift to the level of metaplectic-c bundles, so that  $Q_\sS|_\sC=\sC\times\Mp^c(V_2)$.  The induced one-form $\delta_\sS$ on $Q_{\sS}$ takes the form $$\delta_{\sS}=-\frac{1}{i\hbar}p_c dc+\frac{1}{2}\eta_*\vartheta_0$$ over $\sC$.  Therefore the horizontal lift of $\txi_K$ to $Q_{\sS}|_{\sC}$ is $$\hxi_K(z(c),I)=\lb{\xi_K(z(c)),\overline{\kappa}\oplus\frac{\lambda^2}{i\hbar}},\ \ \forall z(c)\in\sC.$$  Since the $\mathfrak{mp}^c(V_2)$ component is constant, the orbit of $\hxi_K$ through $(z_0,I)\in Q_\sS$ is  $$\widehat{\psi}_K^t(z_0,I)=\lb{\psi_K^t(z_0),\exp\lb{t\overline{\kappa}\oplus\frac{\lambda^2t}{i\hbar}}}.$$  We can write the $\Mp^c(V_2)$ component as $\exp(t\kappa\oplus 0)e^{\lambda^2t/i\hbar}$, where $\exp(t\overline{\kappa}\oplus 0)\in\Mp(V_2)\subset\Mp^c(V_2)$ and $e^{\lambda^2t/i\hbar}\in U(1)\subset\Mp^c(V_2)$.

The parameters of $\exp(t\overline{\kappa}\oplus 0)$ have the form $(\exp{t\overline{\kappa}},\mu(t))$ where $\mu(t)^2\Det_{\C}C_{\exp(t\overline{\kappa})}=1$.  As a $2\times 2$ complex matrix, $$C_{\exp(t\overline{\kappa})}=\four{\cos\lambda t+\frac{i}{2}\lb{\lambda+\frac{1}{\lambda}}\sin\lambda t}{0}{0}{\cos\lambda t+\frac{i}{2}\lb{\lambda+\frac{1}{\lambda}}\sin\lambda t},$$ which has determinant $$\Det_\C C_{\exp(t\overline{\kappa})}=\lb{\cos\lambda t+\frac{i}{2}\lb{\lambda+\frac{1}{\lambda}}\sin\lambda t}^2.$$  This value circles the origin twice as $t$ ranges from $0$ to $2\pi/\lambda$.  Therefore the parameters of $\exp(t\overline{\kappa}\oplus 0)$ are $$\lb{\exp(t\overline\kappa),\lb{\cos\lambda t+\frac{i}{2}\lb{\lambda+\frac{1}{\lambda}}\sin\lambda t}^{-1}},$$ which describes a closed orbit in $\Mp(V_2)$ with period $2\pi/\lambda$.

Thus the orbit in $Q_\sS$ is closed if and only if the $U(1)$ term is also a closed orbit with period $2\pi/\lambda$.  We require $e^{\lambda^2 t/i\hbar}=1$ when $t=2\pi/\lambda$, or equivalently,  $\frac{\lambda^2}{i\hbar}\frac{2\pi}{\lambda}=-2\pi in$ for some $n\in\Z$.  Rearranging and recalling that $\lambda^2=p_c^2=2\E$ results in $\E=\frac{1}{2}n^2\hbar^2$ for some $n\in\Z$.  

We now recall the argument made in Section \ref{subsec:rescale}.  We are only concerned with the strictly positive quantized energies of $K$, because those correspond to quantized energies of the Delaunay Hamiltonian $D$ on $\TpS$.  After dismissing $n<0$ as redundant, this leaves us with  $$\E_n=\frac{1}{2}n^2\hbar^2,\ \ n\in\N.$$  The positive value $\E_n$ is a quantized energy level of $K$ if and only if $f(\E_n)$ is a quantized energy level of $D$, where the diffeomorphism $f$ is given in equation \ref{eq:diffeo}, and the quantized energy levels of $D$ are precisely the negative quantized energy levels of the hydrogen atom.  At last, we find that the quantized energy levels of the hydrogen atom are $$E_n=f(\E_n)=-\frac{m_ek^2}{2n^2\hbar^2},\ \ n\in\N,$$ which is exactly the quantum mechanical prediction for our model of the hydrogen atom.

\end{document}